\documentclass[12pt,reqno]{article}

\usepackage[usenames]{color}
\usepackage{amssymb}
\usepackage{amsmath}
\usepackage{amsthm}
\usepackage{amsfonts}
\usepackage{amscd}
\usepackage{graphicx}

\usepackage[colorlinks=true,
linkcolor=webgreen,
filecolor=webbrown,
citecolor=webgreen]{hyperref}

\definecolor{webgreen}{rgb}{0,.5,0}
\definecolor{webbrown}{rgb}{.6,0,0}

\usepackage{color}
\usepackage{fullpage}
\usepackage{float}

\usepackage{graphics}
\usepackage{latexsym}
\usepackage{epsf}
\usepackage{breakurl}

\newcommand{\seqnum}[1]{\href{https://oeis.org/#1}{\rm \underline{#1}}}
\def\modd#1 #2{#1\ \mbox{\rm (mod}\ #2\mbox{\rm )}}

\usepackage{enumitem}
\usepackage[linesnumbered,lined,boxed,commentsnumbered]{algorithm2e}
\newcommand{\Oh}[1]{\mathcal{O}\left(#1\right)}

\begin{document}


\theoremstyle{plain}
\newtheorem{theorem}{Theorem}
\newtheorem{corollary}[theorem]{Corollary}
\newtheorem{lemma}[theorem]{Lemma}
\newtheorem{proposition}[theorem]{Proposition}

\theoremstyle{definition}
\newtheorem{definition}[theorem]{Definition}
\newtheorem{example}[theorem]{Example}
\newtheorem{conjecture}[theorem]{Conjecture}

\theoremstyle{remark}
\newtheorem{remark}[theorem]{Remark}
\begin{center}
\vskip 1cm{\Large\bf 
On one of Erd\H{o}s' Problems---An Efficient Search for Benelux Pairs
}
\vskip 1cm
\large
Christian Hercher\\
Institut f\"{u}r Mathematik\\
Europa-Universit\"{a}t Flensburg\\
Auf dem Campus 1c\\
24943 Flensburg\\
Germany \\
\href{mailto:christian.hercher@uni-flensburg.de}{\tt christian.hercher@uni-flensburg.de} \\
\end{center}

\vskip .2 in
\begin{abstract}
Erd\H{o}s asked for positive integers $m<n$, such that $m$ and $n$ have the same set of prime factors, $m+1$ and $n+1$ have the same set of prime factors, and $m+2$ and $n+2$ have the same set of prime factors. No such integers are known. If one relaxes the problem and only considers the first two conditions, an infinite series of solutions is known: $m=2^k-2$, $n=(m+1)^2-1=2^k \cdot m$ for all integers $k\geq 2$. One additional solution is also known: $m=75=3\cdot 5^2$ and $n=1215=3^5 \cdot 5$ with $m+1=76=2^2\cdot 19$ and $n+1=1216=2^6 \cdot 19$. No other solutions with $n<2^{32}\approx 4.3\cdot 10^9$ were known.

In this paper, we discuss an efficient algorithm to search for such integers, also known as Benelux pairs, using sieving and hashing techniques. Using highly parallel functioning algorithms on a modern consumer GPU, we could confirm the hitherto known results within a minute of computing time. Additionally, we have expanded the search space by a factor of more than $2^{16}$ and found no further solutions different from the infinite series given above up to $1.4\cdot 10^{12}>2^{40}$. 

For the analogous problem of integers $m<n$ with $m$ and $n+1$ having the same set of prime factors and $m+1$ and $n$having the same set of prime factors, the situation is very similar: An infinite series and one exceptional solution with $n\leq 2^{22}+2^{12}\approx 4.2\cdot 10^6$ were known. We prove that there are no other exceptional solutions with $n<1.4\cdot 10^{12}$.
\end{abstract}

\section{Introduction}
In \cite{Erdos96}, Erd\H{o}s wrote:
\begin{quote}
Balsubramanian called to my attention another old problem of mine which I had forgotten. Can there ever be two distinct integers $m$ and $n$ for which for $0\leq i\leq 2$ the integers $m+i$ and $n+i$ have the same prime factors? The answer is Yes for a more restricted problem. For $m, n$ and $m+1, n+1$ we can take $m=2(2^k-1), n+1=(m+1)^2$. Then $m$ and $n$ also have the same prime factors. Are there other integers $m, n$ and $m+1, n+1$ with the same prime factors?
\end{quote}
{\small \textit{Note: variable names have been changed in this quote for improved readability.}}

Clearly, for all $k\geq 1$,
\begin{itemize}
\item the terms $m=2\cdot (2^k-1)$ and $n=(m+1)^2-1$ are positive integers with $m<n$,
\item $m$ and $n=(m+1)^2-1=m(m+2)=2\cdot (2^k-1) \cdot 2^{k+1}$ have the same prime factors,
\item and $m+1$ and $n+1=(m+1)^2$ have the same prime factors, too. 
\end{itemize}
  Thus, all these pairs are solutions to the restricted problem. (But, since $m+2=2^{k+1}$ is a power of two and $n+2=(2^k-1) \cdot 2^{k+2} + 2 \equiv 2 \pmod{4}$ is not, none of these pairs is a solution to the original problem.)

Section~B19 of \cite{Guy} lists the restricted problem, too. Here, the information is provided that M\c{a}kowski found the exceptional solution $m=3\cdot 5^2$, $n=3^5\cdot 5$ with $m+1=2^2\cdot 19$ and $n+1=2^6\cdot 19$. (But $n+2=7\cdot 11$ and $m+2=1217$ is a prime.) No other exceptional solution is known.

In 2011, the mathematical competition Benelux Olympiad coined in their Problem~1 the term Benelux pair for positive integers $m < n$ when $m$ and $n$, as well as $m+1$ and $n+1$ share the same prime factors.

The original problem from Erd\H{o}s cited above is listed as Number \#850 on Bloom's list~\cite{Bloom850} of Erd\H{o}s' problems. 

A variant of this problem is mentioned in Section~B19 of \cite{Guy}. This variant asks for two positive integers $m<n$ such that now $m$ and $n+1$ should share the same prime factors and $m+1$ and $n$ do as well. The situation for this variant is nearly the same as with Benelux pairs. An infinite series of solutions is given by $m=2^k+1$, $n=m^2-1$. With this setting $m$ and $n+1=m^2$ have the same prime factors, and $m+1=2^k+2=2\cdot(2^{k-1}+1)$ and $n=(m-1)(m+1)=2^k\cdot 2 \cdot(2^{k-1}+1)$ have the same prime factors. There is also one known exceptional solution. According to Guy \cite{Guy}, if $m=5\cdot 7$, $n+1=5^4\cdot 7$, then $n=2\cdot 3^7$, $m+1=2^2\cdot 3^2$. In an analogous way, this paper uses the term \emph{Benelux pair of second kind} for positive integers $m<n$ with $m$ and $n+1$ having the same prime factors, and $m+1$ and $n$ having the same prime factors.

\subsection{Computational verification and contributions from Schott, \linebreak Ehrenstein, Nomoto, and Wasserman}
According to the Online Encyclopedia of Integer Sequences \cite{OEISA343101}, Schott computed all Benelux pairs $(m,n)$ with $n$  up to $2^{30}$. Ehrenstein has expanded the list up to $n<2^{32}$.

For the variant of the problem, Nomoto and Wasserman computed all Benelux pairs of second kind in the OEIS sequence \seqnum{A088966} up to $n\leq 2^{22}+2^{12}$.

\subsection{Contribution of this work}
In this paper, we present an algorithm which computes all Benelux pairs and all Benelux pairs of second kind up to $m<n< S$ with a time complexity of $\Oh{S \log S}$  and a space complexity of $\Oh{S}$. Since the desktop PC used for the computation had a limited amount of RAM, this algorithm is impractical for larger values of~$S$. Therefore, we provide a second algorithm, which only uses a separately changeable amount of space of $\Oh{s}$. But in a time-space-trade-off it needs a time of $\Oh{\frac{1}{s} \cdot S^2}$ to compute, assuming $s\in \Omega(\log S \log\log S)$. This second algorithm has the advantage that the parts that need the most computations can be highly parallelized. Thus, we used the GPU of the PC in hand to extend the search for such pairs by a large amount.

\subsection{Outline}
In Section~\ref{sec:Algos}, we construct and discuss in detail the two aforementioned algorithms before we conclude in Section~\ref{sec:Results} with a discussion of further improvements for high parallel computing, and provide the results of our computation.
   
\section{Algorithms}\label{sec:Algos}
\subsection{Finding Benelux pairs quickly}
Both problems---finding Benelux pairs of first and second kind---consider pairs of positive integers and their respective prime factors (ignoring multiplicities). Thus, a good way to start would be to calculate for every positive integer~$n\leq S$ the set of primes which divide~$n$. A first inside in these problems is that we do not have to know which primes exactly divide $n$, but whether two such sets are identical or not. Hence, we do not need to compute the whole prime factorization of~$n$ if we only compute a number $r(n)$ with the property that $r(n)=r(m)$, if and only if $n$ and $m$ have the same prime factors. Such a function $r$ is the \emph{radical} of~$n$: $rad(n):=\prod_{p\mid n} p$. And this function has the great advantage that if we want to determine its value for all elements in a large interval, it can be computed efficiently via sieving.

This can be achieved using Algorithm~\ref{Algo_Sieving_Radicals}: First, for every integer~$n$ in the interval $[start,start+length[$ we initialize the value $r[n]$ as $r[n]:=n$. Then we will subsequently divide this value by all primes~$p$ which have a multiplicity of $>1$ in the prime factorization of~$n$---and do this multiple times until the prime~$p$ remains with multiplicity~1 in the resulting quotient. To do this, we handle each prime~$p$ one at a time. The largest prime we have to consider is at maximum $\sqrt{start+length-1}$ since for all larger primes~$p$, the integer $p^2$ cannot be a divisor of any number in the considered interval. For all primes $p$ up to this limit, we identify all integers in the considered interval that are divisible by $p^2$, $p^3$, $\dots$, and divide the value $r[n]$ by $p$. Thus, if the prime~$p$ has a multiplicity of $m\geq 2$ in the prime factorization of an integer~$n$, it is $n$ divisible by $p^2$, $p^3$, $\dots$, $p^m$, but not by $p^{m+1}$. So, in this process, we divide the value $r[n]$ exactly $m-1$ times by $p$, hence after this process $p\mid r[n]$ but $p^2\nmid r[n]$. When we have done this for all primes~$p$, we have to consider, it is $r[n]=rad(n)$ for all~$n$ in the interval. Thus, Algorithm~\ref{Algo_Sieving_Radicals} in fact computes the radicals of all integers in the given interval.

\begin{center}
\begin{algorithm}[H]
 \SetAlgoLined
 \SetKw{Kwcont}{continue}
  \KwIn{$start$ and $length$ of an interval to screen}
  \KwOut{for every $n$ in $[start; start+length[$ the radical $r[n]=rad(n)$}
  \ForAll{$n \in [start, start+length[$}
  {Set $r[n] \leftarrow n$\;}
  \BlankLine
  \ForAll{primes $p\leq \sqrt{start+length-1}$}
  {
      Set $exponent \leftarrow 2$\;
      \While{$p^{exponent}\leq start+length-1$ }
      {
	  	Set $res \leftarrow start \mod{p^{exponent}}$ \;
	  	\If{$res = 0$}
	          {
	               Set $res \leftarrow p^{exponent}$\;
	           }
	           Set $shift \leftarrow p^{exponent}-res$\;
	           \tcp{So $start+shift$ is the first number in the interval which is divisible by $p^{exponent}.$}
	           
	           \If(\tcp*[h]{nothing more to do for this prime}){$shift>length$}{\Kwcont \tcp{Move to next prime.}}
	           \BlankLine
	           
	           \ForAll{integers $0\leq k\leq \frac{length-1-shift}{p^{exponent}}$}
	           {
	               Divide $r[start+shift+k\cdot p^{exponent}]$ by $p$\;
	           }
	           
	           Increment $exponent$\;
      }
  }
    \Return{the vector $r$}\;
\caption{Computing radicals of all integers in an interval via sieving.}
\label{Algo_Sieving_Radicals}
\end{algorithm}
\end{center}

For every prime~$p$, Algorithm~\ref{Algo_Sieving_Radicals} has to do $\Oh{\frac{length}{p^2}+\frac{length}{p^3}+\dots}=\Oh{length \cdot \frac{1}{p^2}}$ dividing steps. Since $\sum_p \frac{1}{p^2}<\sum_{n=1}^{\infty} \frac{1}{n}=\zeta(2)$ converges, there are only $\Oh{length}$ divisions needed. Furthermore, let $S:= start+length
1$. Then for each prime~$p$, $\Oh{\frac{\log S}{\log p}}$ elementary update operations are requiered. Since $\sum_{p\leq \sqrt{S}} \frac{1}{\log p}<\int_{x=2}^{\sqrt{S}} \frac{\text{d}x}{\log x}=\log\log \sqrt{S} - \log\log 2 < \log\log S$ that are, at most, $\Oh{\log S \cdot \log\log S}$ such $\Oh{1}$ update steps in total. Thus, if the set of primes $p\leq \sqrt{S}$ is precomputed,  Algorithm~\ref{Algo_Sieving_Radicals} runs for, at most, $\Oh{length+\log S \cdot \log\log S}$ steps and needs memory of size $\Oh{length}$, since obviously $rad(n)$ is held in memory for every integer in the interval. If $start=1$, hence $length=S$, this reduces to a time complexity of $\Oh{S}$ and a space complexity of $\Oh{S}$ as well.
 
The next insight we observe is that both problems do not ask for the radicals of integers $m<n$ specifically, but only, if the sets $\{rad(m),rad(m+1)\}$ and $\{rad(n),rad(n+1)\}$ are equal: In the case of $rad(m)=rad(n)$ and $rad(m+1)=rad(n+1)$, we have found a Benelux pair $(m,n)$; and if $rad(m)=rad(n+1)$ and $rad(m+1)=rad(n)$, $(m,n)$ is a Benelux pair of second kind. Thus, after computing the radicals for every $n<S$, we build the sets $S_n:=\{rad(n),rad(n+1)\}$.

A simple way to find duplicates in a list is to sort the list and then traverse the sorted list. To do this, we define an easily computable order relation between sets with two elements: Say $\{a_1,a_2\} < \{b_1,b_2\}$ with $a_1<a_2$ and $b_1<b_2$, if and only if $a_1<b_1$ or ($a_1=b_1$ and $a_2<b_2$). Then, for two such sets $A$ and $B$, we clearly always have exactly one of the three possibilities $A<B$, $B<A$ or $A=B$. With this order relation, we can sort the sets $S_n$ and find collisions, in other words, integers $m<n$ with $S_m=S_n$, and hence Benelux pairs of first or second kind. Since every Benelux pair of either kind gives rise to a collision, this approach finds every such pair in the given search space.

The whole procedure is summarized in Algorithm~\ref{Algo_Find_Benelux_fast}.

\begin{center}
\begin{algorithm}[H]
 \SetAlgoLined
  \SetKw{Kwprint}{print}
  \KwIn{upper limit $S$}
  \KwOut{all Benelux pairs $(m,n)$ of first and second kind with $m<n<S$}
  \ForAll(\tcp*[h]{Use Algorithm~\ref{Algo_Sieving_Radicals} for this loop.}){$1\leq n \leq S$}
  {Compute $r[n]$\;}
  
  \BlankLine
  \ForAll{$1\leq n < S$}
  {
      Set $S_n \leftarrow \{r[n], r[n+1]\}$\;
      Set $list[n] \leftarrow (n,S_n)$\;
  }
  \BlankLine
  
  Sort list with respect to $(m,S_m) < (n, S_n) \iff S_m<S_n$\;
  
  \ForAll{$1\leq i < S-1$}
  {
  	Set $(m, S_m) \leftarrow list[i]$\;
  	Set $(n, S_n) \leftarrow list[i+1]$\;
  	\If(\tcp{duplicate found}){$S_m=S_n$}
  	{
  		\eIf{$r[m]=r[n]$}
  		{\Kwprint{Benelux pair of 1st kind found: $(m,n)$.}}
  		{\Kwprint{Benelux pair of 2nd kind found: $(m,n)$.}}
  	}
  }
  
\caption{Finding all Benelux pairs of first and second kind up to $S$.}
\label{Algo_Find_Benelux_fast}
\end{algorithm}
\end{center}

Computing all radicals with Algorithm~\ref{Algo_Sieving_Radicals} requires $\Oh{S}$ steps. Building the sets $S_n$ and the list can be achieved in $\Oh{S}$, too. Sorting this list takes $\Oh{S\log S}$ steps and the search for duplicates in this sorted list is completed in time $\Oh{S}$. Combining this gives a time complexity of $\Oh{S \log S}$. And since a constant number of arrays of size $S$ are needed, Algorithm~\ref{Algo_Find_Benelux_fast} needs a memory amount of $\Oh{S}$.

So, if memory is not a limiting factor, the search for Benelux pairs can be done asymptotically as fast as sorting.

\subsection{Breaking the space limitation barrier}
Unfortunately, memory quickly becomes a limiting factor: The size of RAM in modern desktop computers is in the order of gigabytes. So, if the search should be extended beyond $2^{32}$, this is no longer sufficient. Thus, we have to come up with new ideas to get an effective search algorithm even for a larger search spaces.

If we do not have sufficient memory to hold all results of our computations, we have no other choice but to divide the search space into smaller chunks. Let~$s$ be the size of such a chunk; in particular, let the chunk $C_i$ be the set of all integers in the interval $[1+i\cdot (s-1), 1+(i+1) \cdot (s-1)]$. (Two successive chunks have to have an overlap of one number since we also need the value of $rad(n+1)$ for every integer~$n$. This value cannot be computed for the last number in a chunk---so the last number has to be considered at the start of the next chunk a second time.) 

We can compute the radicals with Algorithm~\ref{Algo_Sieving_Radicals} and build the sets $S_n$ for all integers in this chunk in the same way as above in Algorithm~\ref{Algo_Find_Benelux_fast}. We also could find duplicate sets within a chunk in the same way as above. But how can pairs $(m,n)$ with $S_m=S_n$ be found if $m$ and $n$ are in different chunks? For this problem (and in fact for finding duplicates within a chunk, too) we use another approach.

Another efficient way of searching for duplicates is to use hashing: Insert the hash values of the ojects under consideration into a hash table and whenever one gets a collision one checks whether the objects in question where equal in the first place. If one uses a hash function and a table size which leads to few collisions for different objects, this should be a fast alternative. This gives us an approach, as used in Algorithm~\ref{Algo_Find_Benelux_with_Hash}, for finding all Benelux pairs $(m,n)$ of both kinds with $m<n$ and $n$ being in the currently examined chunk $C_i$: Algorithm~\ref{Algo_Find_Benelux_with_Hash} computes the radical $r[n]$ and set $S_n$ for all numbers~$n$ in the given chunk $C_i$ as above. But then it computes hash values $h_n$ for every such set $S_n$ and uses them to insert these data into a previously empty hash table. If there is a collision, the two sets are tested to see if they are equal. (If they are, then a Benelux pair has been found.) Open addressing can be used to find another slot to insert this data point. Now fix this hash table and recompute the radicals, sets, and hash values for all previous chunks $C_j$ with $j<i$. For every number~$m$ in such a previous chunk, search whether $h_m$ is in the above generated hash table. If so, test if the corresponding sets are equal. 

\begin{center}
\begin{algorithm}[H]
 \SetAlgoLined
  \SetKw{Kwprint}{print}
  \KwIn{the chunk $C_i$ of integers}
  \KwOut{all Benelux pairs $(m,n)$ of first and second kind with $m<n$ and $n\in C_i$}
  
  Reset the hash table to an empty one\;
  \BlankLine

  \ForAll(\tcp*[h]{Use Algorithm~\ref{Algo_Sieving_Radicals} for this loop.}){$n \in C_i$}
  {Compute $r[n]$\;}
  
  \BlankLine
  \ForAll{$n \in C_i$}
  {
      Set $S_n \leftarrow \{r[n], r[n+1]\}$\;
      Compute a hash value $h_n$ from $S_n$\;
      Insert $(n,h_n)$ with key value $h_n$ in the hash table\;
      \If{there is a collision while $(n,h_n)$ is inserted \label{Algo4:insertHashTable}}
      {
      	Let $(m,h_m)$ be the date already in the hash table\;
      	\If{$S_m=S_n$}{Test whether $(m,n)$ is a Benelux pair of 1st or 2nd kind and print it\;}
      	Use open addressing to find an unoccupied slot for the current  data\;
      }
  }
  \tcp{This finds all Benelux pairs $(m,n)$ of both kinds with $m,n\in C_i$.}
  \BlankLine
   
  \ForAll{$0\leq j < i$ \label{Algo4:otherChunks}}
  {
  	\ForAll(\tcp*[h]{Use Algorithm~\ref{Algo_Sieving_Radicals} for this loop.}){$m \in C_j$}
  	{Compute $r[m]$\;}
  	\ForAll{$m \in C_j$}
	 {
	      Set $S_m \leftarrow \{r[m], r[m+1]\}$\;
	      Compute a hash value $h_m$ from $S_m$\;
	      Search for the key value $h_m$ in the hash table\;
	      \If{the search is successful \label{Algo4:searchInHashTable}}
	      {
	      	Let $(n,h_n)$ be the date in the hash table\;
	      	\If{$S_m=S_n$}{Test whether $(m,n)$ is a Benelux pair of 1st or 2nd kind and print it\;}
	      	Test until an unoccupied slot is found \tcp*{Then there can't be another entry with same hash value.}
	      }
	}
	\tcp{This finds all Benelux pairs $(m,n)$ of both kinds with $m\in C_j$ and $n\in C_i$.}
  }
  
\caption{Finding all Benelux pairs $(m,n)$ of both kinds with $m<n$ and $n\in C_i$.}
\label{Algo_Find_Benelux_with_Hash}
\end{algorithm}
\end{center}

If $(m,n)$ is a Benelux pair of either kind with $m$ and $n$ from the same chunk, Algorithm~\ref{Algo_Find_Benelux_with_Hash} will identify it in Step~\ref{Algo4:insertHashTable} where the collision is detected while $(n,h_n)$ is being inserted into the hash table. And if $m$ is from a chunk $C_j$ with $j<i$, this pair will be found in  Step~\ref{Algo4:searchInHashTable} while searching for the hash value~$h_m$ of $S_m$ in the hash table. Thus, all Benelux pairs $(m,n)$ with $m<n$ and $n\in C_i$ are found. Since the only values stored in the hash table are hash values $h_n$ from numbers $n$ in the currently investigated chunk $C_i$, no other (potential) solutions than the stated ones are reported. Hence, the algorithm is correct.

As the length of each chunk is~$s$, computing the radicals of all numbers in a chunk using Algorithm~\ref{Algo_Sieving_Radicals} needs $\Oh{s + \log S\log\log S}$ steps. Assuming that $s$ is in $\Omega(\log S\log\log S)$ this reduces to $\Oh{s}$. The sets and hash values of all integers in a given chunk can be computed in $\Oh{s}$ steps. Let the hash function used be constructed so that two randomly chosen different sets in a family of $s$ such sets have the same hash value only with the probability~$p<1$. If one uses a hash table with size $t\cdot s$, the probability of collisions for different hash values is $\frac{1}{t}$. Thus, while inserting a data point in the hash table in Step~\ref{Algo4:insertHashTable} of the algorithm, the expected number of slots to be tested until an empty one is found, where the new data can be stored, is $\Oh{\frac{1}{p} \cdot \frac{t}{t-1}}$. So, with constants $p$ and $t$, the hash table can be filled with the data for all numbers in one chunk in $\Oh{s}$ steps. 

Now the loop beginning in Step~\ref{Algo4:otherChunks} of Algorithm~\ref{Algo_Find_Benelux_with_Hash} has to be computed for all $\Oh{\frac{S}{s}}$ chunks which came before the current one. In every iteration of this loop, the radicals, sets, and hash values for all numbers in a chunk $C_j$ can be computed in $\Oh{s}$ steps, as above. In the same way as in Step~\ref{Algo4:insertHashTable}, searching for a duplicate to a single data point can be done in an expected running time of $\Oh{1}$, which is $\Oh{s}$ for all elements of a chunk $C_j$. Hence, all computations for such a chunk $C_j$ need $\Oh{s}$ steps and the whole loop beginning in Step~\ref{Algo4:otherChunks} has a running time of $\Oh{S}$. Since $s\leq S$, Algorithm~\ref{Algo_Find_Benelux_with_Hash} needs $\Oh{S}$ steps.

But to find all Benelux pairs $(m,n)$ of both kinds with $m<n$, one has to run Algorithm~\ref{Algo_Find_Benelux_with_Hash} for all $\frac{S}{s-1}$ chunks consisting of integers $\leq S$. Thus, one can identify all Benelux pairs of both kinds in $\Oh{\frac{S^2}{s}}$ steps. And since at any one time only a constant number of lists of size~$s$ is needed, these computations can be achieved wile only using memory of size~$\Oh{s}$.

\section{Implementation and results}\label{sec:Results}
Already the first steps in implementing Algorithm~\ref{Algo_Find_Benelux_fast} showed that the memory limitations of the desktop computer available to the author would quickly limit the search space. Thus, the most effort was invested in implementing and optimizing Algorithm~\ref{Algo_Find_Benelux_with_Hash}. To get the best possible performance, one should chose the chunk size~$s$ to be as large as possible. The author went for $2^{27}$. This allowed for a problem-free use of a large part of the memory available. Now the most computationally intensive parts of the algorithm have the advantage of being highly parallelizable. Thus, we optimized the code to use this and utilize the significant parallelism that is possible on a modern GPU.

\subsection{Optimization for parallel computing}
A number of optimizations were applied while implementing Algorithm~\ref{Algo_Find_Benelux_with_Hash}. Some decisive ones---in the sense of accelerating the code--- are listed below:

\begin{itemize}
\item A key part in the algorithm is the computation of the radicals. Given a chunk $C_i$, one has to identify all multiples of $p^k$ and divide the respective values by $p$. If the offset (which gives the first number in the chunk that is divisible by $p^k$) is computed for a given prime power $p^k$, then these divisions can be done all in parallel since they are independent of each other. Furthermore, if $p<q$ are two primes with $p^4\geq start+length$, then no number in the given chunk can be divisible by $p^2$ and $q^2$ at the same time. Thus, no two such primes can influence the computation for the same number and so with all of these primes the sieving process can be done in parallel. (If $p^3\geq start+length$, we at least know that no number in the chunk is divisible by $p^3$ or higher powers, thus we only have to check for multiples of $p^2$.)
\item For the prime $p=2$, we can use a built-in function provided by CUDA (a programming language and C++ extension to write code for NVIDIA GPUs): $ctz$, that is ``count trailing zeros''. This function returns the number of zeros that a binary integer ends in, hence $ctz(n)=x$ is equivalent to $2^x\mid n$, but $2^{x+1}\nmid n$. Thus, for even~$n$, we have to divide $n$ by $2^{x-1}$ to only let one factor~2 remain. (For even $n$ that are not a multiple of~4, there is nothing to do in this context.) So, for all $n$ which are multiples of~4,  we compute $ctz(n)$ and then shift $n$ by $ctz(n)-1$ bits to the right. Both operations are very fast and can be performed in parallel for all~$n$. Then the prime $p=2$ is fully processed. 
\item After the radicals of all numbers in the given chunk are computed, we have to compute the hash values of the respective sets. For this we use a commutative hash function, which sends a pair $(rad(n), rad(n+1))$ of 64 bit integers to an integer $0\leq h_n < 2^{29}=4s$. As this is the size of the hash table, too, the initial hash adresses for two different sets should be equal only in $\frac{1}{4}$ of the cases. Thus, since we use open addressing, on average only a small number of slots has to be examined to find an unoccupied one to insert a new data point. In parallel, we can examine the slot at the initial address, the next one, and so on, writing in shared memory whichever address is the first with an unoccupied slot. (In our tests, the best performance was achieved by examine 4 addresses in parallel.) In the same way as the insertion of new data into the hash table the search for duplicates can be handled. And, of course, calculating the hash values, inserting data into the hash table and searching in the hash table can be done in parallel for all elements of a chunk. One only has to ensure that no two values are inserted at the same time at the same slot in the hash table. This is be done with built-in thread-locking mechanisms.
\item Since on the GPU only 32 bit integers are supported natively, we had to implement our own 64 bit arithmetic. By doing so we could optimize for our special circumstances. For example, in the recurring task of reducing a 64 bit integer $a_1\cdot 2^{32}+a_0$ modulo a prime $p<2^{32}$, we computed $((a_1 \mod p) \cdot (2^{32} \mod p) + (a_0 \mod p)) \mod p$, where the factor $(2^{32} \mod p)$ could be precomputed in advance and only has then to be read from memory.
\item Since communication between host memory and device memory (in other words, RAM vs.\ memory on the GPU) could result in a bottleneck in such computations, we reduced the communication to a bare minimum: After sending some initial data---such as the list of primes to be used in the sieving process---to the GPU, transfer of data was limited to primarily the chunk number. All other data---such as which numbers are in this chunk and which of them are part of a Benelux pair---are computed entirely on the GPU. Possible findings of Benelux pairs are written into a buffer and sent back to the host system, which writes them into a file. 
\end{itemize}

\subsection{List of findings}
Using the optimized implementation of Algorithm~\ref{Algo_Find_Benelux_with_Hash}, we were able to extend the previously considered search space by a large factor: In approximately one minute of computing time on our standard desktop PC with a consumer graphics card from 2021, we reached $S=2^{32}$, and thus the upper limit of previous searches. Since the algorithm for fixed~$s$ is quadratic in the upper search limit~$S$, it took approximately one and a half months to reach $S=2^{40}\approx 1.0995\cdot 10^{12}$. In fact, we searched a little further. From this we can conclude:

\begin{theorem}
Let $m<n<1.4 \cdot 10^{12}$ be positive integers. Then
\begin{enumerate}[label=\alph*)]
	\item $(m,n)$ is a Benelux pair (of first kind) if and only if
	\begin{itemize}
		\item $m=75$ and $n=1215$ or
		\item $m=2^k-2$ and $n=2^{2k}-2^{k+1}$ for an integer $2\leq k\leq 20$ and
	\end{itemize}
	\item $(m,n)$ is a Benelux pair of second kind if and only if
	\begin{itemize}
		\item $m=35$ and $n=4374$ or
		\item $m=2^k+1$ and $n=2^{2k}+2^{k+1}$ for an integer $0\leq k\leq 20$.
	\end{itemize}
\end{enumerate}

Table~\ref{Tab:BeneluxPairs} gives all Benelux pairs with their prime factorizations and Table~\ref{Tab:BeneluxPairs2nd} does the same for Belenux pairs of second kind. 
\end{theorem} 

\begin{table}[H]
\footnotesize
\begin{center}
\begin{tabular}{r|r||r|r||l} 
$m$ & $n$ & $m+1$ & $n+1$ & Remarks\\
\hline
\hline
$75=$                 & $1215=$                        & $76=$                                                      & $1216=$                    & exceptional\\
$3\cdot 5$          & $3^5 \cdot 5$                & $2^2\cdot 19$                                         & $2^6 \cdot 19$         & \\
\hline 
$2=$                   & $8=$                               & $3=$                                                        & $9=$                         & $k=\phantom{1}2$\\
$2$                      & $2^3$                            & $3$                                                           & $3^2$                       &  \\ 
\hline
$6=$                    & $48=$                            & $7=$                                                        & $49=$                       & $k=\phantom{1}3$\\
$2\cdot 3$           & $2^3 \cdot 3$                & $7$                                                           & $7^2$                       & \\
\hline
$14=$                  & $224=$                          & $15=$                                                      & $225=$                      & $k=\phantom{1}4$\\
$2\cdot 7$           & $2^5 \cdot 7$                & $3\cdot 5$                                               & $3^2 \cdot 5^2$       & \\
\hline
$30=$                  & $960=$                          & $31=$                                                      & $961=$                      & $k=\phantom{1}5$\\
$2\cdot 3 \cdot 5$& $2^6 \cdot 3 \cdot 5$  & $31$                                                         & $31^2$                     & \\
\hline
$62=$                  & $3968=$                        & $63=$                                                     & $3969=$                    & $k=\phantom{1}6$\\
$2\cdot 31$         & $2^7 \cdot 31$              & $3^2\cdot 7$                                          & $3^4 \cdot 7^2$       & \\
\hline
$126=$                & $16118=$                      & $127=$                                                     & $16119=$                & $k=\phantom{1}7$\\
$2\cdot 3^2 \cdot 7$ & $2^8 \cdot 3^2 \cdot 7$ & $127$                                              & $127^2$                   & \\
\hline
$254=$                & $65024=$                      & $255=$                                                     & $65025=$               & $k=\phantom{1}8$\\
$2\cdot 127$       & $2^9 \cdot 127$            & $3 \cdot 5 \cdot 17$                                 & $3^2 \cdot 5^2 \cdot 17^2$  & \\
\hline
$510=$                & $261120=$                    & $511=$                                                     & $261121=$             & $k=\phantom{1}9$\\
$2 \cdot 3 \cdot 5 \cdot 17$ & $2^{10} \cdot 3 \cdot 5 \cdot 17$  & $7\cdot 73$          & $7^2\cdot 73^2$    & \\
\hline
$1022=$              & $1046528=$                  & $1023=$                                                   & $1046529=$           & $k=10$\\
$2 \cdot 7\cdot 73$ & $2^{11} \cdot 7\cdot 73$ & $3\cdot 11 \cdot 31$                       & $3^2\cdot 11^2 \cdot 31^2$    & \\
\hline
$2046=$              & $4190208=$                  & $2047=$                                                   & $4190209=$           & $k=11$\\
$2 \cdot3\cdot 11 \cdot 31$ & $2^{12} \cdot 3\cdot 11 \cdot 31$ & $23 \cdot 89$       & $23^2\cdot 89^2$  & \\
\hline
$4094=$              & $16769024=$               & $4095=$                                                   & $16769025=$         & $k=12$\\
$2 \cdot 23\cdot 89$ & $2^{13} \cdot 23\cdot 89$ & $3^2 \cdot 5 \cdot 7 \cdot 13$    & $3^4 \cdot 5^2 \cdot 7^2 \cdot 13^2$  & \\
\hline
$8190=$              & $67092480=$               & $8191=$                                                   & $67092481=$         & $k=13$\\
$2 \cdot 3^2 \cdot 5 \cdot 7 \cdot 13$ & $2^{14} \cdot 3^2 \cdot 5 \cdot 7 \cdot 13$ & $8191$ & $8191^2$  & \\
\hline
$16382=$            & $268402689=$             & $16383=$                                                 & $268402689=$       & $k=14$\\
$2 \cdot 8191$    & $2^{15} \cdot 8191$   & $3 \cdot 43 \cdot 127$                             & $3^2 \cdot 43^2 \cdot 127^2$  & \\
\hline
$32766=$            & $1073676288=$            & $32767=$                                                 & $1073676289=$    & $k=15$\\
$2 \cdot 3 \cdot 43 \cdot 127$ & $2^{16} \cdot 3 \cdot 43 \cdot 127$ & $7 \cdot 31 \cdot 151$ & $7^2 \cdot 31^2 \cdot 151^2$  & \\
\hline
$65534=$            & $4294836224=$            & $65535=$                                                 & $4294836225=$    & $k=16$\\
$2 \cdot 7 \cdot 31 \cdot 151$ & $2^{17} \cdot 7 \cdot 31 \cdot 151$ & $3 \cdot 5 \cdot 17 \cdot 257$ & $3^2 \cdot 5^2 \cdot 17^2 \cdot 257^2$  & \\
\hline
$131070=$          & $17179607040=$         & $131071=$                                               & $17179607041=$   & $k=17$\\
$2 \cdot 3 \cdot 5 \cdot 17 \cdot 257$ & $2^{18} \cdot 3 \cdot 5 \cdot 17 \cdot 257$ & $131071$ & $131071^2$  & \\
\hline
$262142=$          & $68718952448=$         & $262143=$                                               & $68718952449=$   & $k=18$\\
$2 \cdot131071$ & $2^{19} \cdot 131071$& $3^3 \cdot 7 \cdot 19 \cdot 73$ & $3^6 \cdot 7^2 \cdot 19^2 \cdot 73^2$ &\\
\hline
$524286=$         & $274876858368=$         & $524287=$                                               & $274876858369=$   & $k=19$\\
$2\cdot 3^3 \cdot 7 \cdot 19 \cdot 73$ & $2^{20} \cdot 3^3 \cdot 7 \cdot 19 \cdot 73$ & $524287$                                                & $524287^2$              & \\
\hline
$1048574=$        & $1099509530624=$      & $1048575=$                                            & $1099509530625=$ & $k=20$\\
$2\cdot 524287$ & $2^{21} \cdot 524287$ & $3\cdot5^2\cdot 11 \cdot 31 \cdot 41$ & $3^2\cdot5^4\cdot 11^2 \cdot 31^2 \cdot 41^2$ & \\
\hline
\end{tabular}
\caption{All Benelux pairs $(m,n)$ with $m<n<1.4\cdot 10^{12}$ and their prime factorizations}
\label{Tab:BeneluxPairs}
\end{center}
\end{table}

\begin{table}[H]
\footnotesize
\begin{center}
\begin{tabular}{r|r||r|r||l}
$m$ & $n+1$ & $m+1$ & $n$ & Remarks\\
\hline
\hline
$35=$                 & $4375=$                        & $36=$                                                      & $4374=$                    & exceptional\\
$5\cdot 7$          & $5^4 \cdot 7$                & $2^2\cdot 3^2$                                      & $2 \cdot 3^7$           & \\
\hline 
$2=$                   & $4=$                               & $3=$                                                        & $3=$                         & $k=\phantom{1}0$\\
$2$                      & $2^2$                            & $3$                                                           & $3$                          &  \\ 
\hline
$3=$                   & $9=$                              & $4=$                                                        & $8=$                         & $k=\phantom{1}1$\\
$3$                      & $3^2$                            & $2^2$                                                      & $2^3$                       &  \\ 
\hline
$5=$                   & $25=$                            & $6=$                                                        & $24=$                       & $k=\phantom{1}2$\\
$5$                      & $5^2$                            & $2 \cdot 3$                                              & $2^3 \cdot 3$           &  \\ 
\hline
$9=$                    & $81=$                            & $10=$                                                     & $80=$                       & $k=\phantom{1}3$\\
$3^2$                 & $3^4 $                            & $2\cdot 5$                                              & $2^4 \cdot 5$           & \\
\hline
$17=$                  & $289=$                          & $18=$                                                     & $288=$                      & $k=\phantom{1}4$\\
$17$                    & $17^2$                           & $2\cdot 3^2$                                         & $2^5 \cdot 3^2$       & \\
\hline
$33=$                  & $1089=$                        & $34=$                                                    & $1088=$                      & $k=\phantom{1}5$\\
$3 \cdot 11$        & $3^2 \cdot 11^2$         & $2 \cdot 17$                                          & $2^6 \cdot 17$           & \\
\hline
$65=$                  & $4225=$                        & $66=$                                                   & $4224=$                    & $k=\phantom{1}6$\\
$5\cdot 13$         & $5^2 \cdot 13^2$         & $2\cdot 3 \cdot 11$                               & $2^7 \cdot 3 \cdot 11$ & \\
\hline
$129=$                & $16641=$                      & $130=$                                                 & $16640=$                & $k=\phantom{1}7$\\
$3\cdot 43$         & $3^2 \cdot 43^2$         & $2 \cdot 5 \cdot 13$                             & $2^7 \cdot 5 \cdot 13$   & \\
\hline
$257=$                & $66049=$                      & $258=$                                                & $66048=$               & $k=\phantom{1}8$\\
$257$                  & $2577^2$                      & $2 \cdot 3 \cdot 43$                            & $2^8 \cdot 3 \cdot 43$  & \\
\hline
$513=$                & $263169=$                    & $514=$                                               & $263168=$             & $k=\phantom{1}9$\\
$3^3 \cdot 19$    & $3^6 \cdot 19^2$        & $2 \cdot 257$                                     & $2^9 \cdot 257$    & \\
\hline
$1025=$              & $1050625=$                  & $1026=$                                                   & $1050624=$           & $k=10$\\
$5^2\cdot 41$     & $5^4\cdot 41^2$         & $2 \cdot 3^3 \cdot 19$                       & $2^{10} \cdot 3^3 \cdot 19$    & \\
\hline
$2049=$              & $4198401=$                  & $2050=$                                                   & $4198400=$           & $k=11$\\
$3\cdot 683$       & $3^2\cdot 683^2$        & $2 \cdot 5^2\cdot 41$                       & $2^{11} \cdot 5^2\cdot 41$  & \\
\hline
$4097=$              & $16785409=$               & $4098=$                                                   & $16785408=$         & $k=12$\\
$17\cdot 241$     & $17^2\cdot 241$          & $2 \cdot 3\cdot 683$                          & $2^{12} \cdot 3\cdot 683$  & \\
\hline
$8193=$              & $67125249=$               & $8194=$                                                   & $67125248=$         & $k=13$\\
$3\cdot 2731$     & $3^2 \cdot 2731^2$    & $2 \cdot 17\cdot 241$                        & $2^{13} \cdot 17\cdot 241$  & \\
\hline
$16385=$            & $268468225=$             & $16386=$                                                 & $268468224=$       & $k=14$\\
$5 \cdot 29 \cdot 113$ & $5^2 \cdot 29^2 \cdot 113^2$ & $2 \cdot 3\cdot 2731$ & $2^{14} \cdot 3\cdot 2731$  & \\
\hline
$32769=$            & $1073807361=$            & $32770=$                                                 & $1073807360=$    & $k=15$\\
$3^2 \cdot 11 \cdot 331$ & $3^4 \cdot 11^2 \cdot 331^2$ & $2 \cdot 5 \cdot 29 \cdot 113$ & $2^{15} \cdot 5 \cdot 29 \cdot 113$  & \\
\hline
$65537=$            & $4295098369=$            & $65538=$                                                & $4295098368=$    & $k=16$\\
$65537$               & $65537^2$                    & $2 \cdot 3^2 \cdot 11 \cdot 331$  &$2^{16} \cdot 3^2 \cdot 11 \cdot 331$& \\
\hline
$131073=$          & $17180131329=$         & $131074=$                                               & $17180131328=$   & $k=17$\\
$3 \cdot 43691$  & $3^2 \cdot 43691^2$  & $2 \cdot 65537$                               & $2^{17} \cdot 65537$  & \\
\hline
$262145=$          & $68720001025=$         & $262146=$                                               & $68720001024=$   & $k=18$\\
$5\cdot 13 \cdot 37 \cdot 109$ & $5^2\cdot 13^2 \cdot 37^2 \cdot 109^2$& $2 \cdot 3 \cdot 43691$ & $2^{18} \cdot 3 \cdot 43691$ &\\
\hline
$524289=$         & $274878955521=$         & $524290=$                                               & $274878955520=$   & $k=19$\\
$3\cdot 174763$ & $3^2 \cdot 174763^2$ & $2 \cdot 5\cdot 13 \cdot 37 \cdot 109$   & $2^{19} \cdot 5\cdot 13 \cdot 37 \cdot 109$              & \\
\hline
$1048577=$        & $1099513724929=$      & $1048578=$                                            & $1099513724928=$ & $k=20$\\
$17 \cdot 61681$ & $17^2 \cdot 61681^2$ & $2 \cdot 3\cdot 174763$             & $2^{20} \cdot 3\cdot 174763$ & \\
\hline
\end{tabular}
\caption{All Benelux pairs $(m,n)$ of second kind with $m<n<1.4\cdot 10^{12}$ and their prime factorizations}
\label{Tab:BeneluxPairs2nd}
\end{center}
\end{table}

\bigskip
\hrule
\bigskip

\noindent 2020 {\it Mathematics Subject Classification}:
Primary 11-04; Secondary 11Y55.

\noindent \emph{Keywords: }Erd\H{o}s problem, Benelux pair, integers with same prime factors. 
\bigskip
\hrule
\bigskip

\noindent (Concerned with sequences \seqnum{A343101} and \seqnum{A088966} of the On-Line Encyclopedia of Integer Sequences.)

\bigskip
\hrule
\bigskip

\end{document}